\documentclass[10pt]{article}

\usepackage{theorem,amssymb,amsmath}

\usepackage{graphicx}

\usepackage{hyperref}

\usepackage{color}                    %%%%%%%%
\definecolor{red}{rgb}{1,0,.2}        %%%%%%%%
\definecolor{cjp}{rgb}{.1,.7,.2}        %%%%%%%%
\definecolor{fmdc}{rgb}{1,0,.8}        %%%%%%%%
\newcommand{\Sh}{{\rm S}}
\newcommand{\Lan}{{\mathcal{L}}}
\newcommand{\T}{{\cal T}}
\newcommand{\xF}{x_{\rm \scriptstyle F}}
\newcommand{\N}{\mathbb{N}}

\advance \topmargin by -\headheight
\advance \topmargin by -\headsep
\evensidemargin \oddsidemargin
\author{J.-P. Allouche \\
CNRS, IMJ-PRG, Sorbonne Universit\'e \\
4 Place Jussieu \\
F-75252 Paris Cedex 05 France \\
{\tt jean-paul.allouche@imj-prg.fr}
\and
F. M. Dekking \\
Delft University of Technology \\
Faculty EEMCS, P.O.~Box 5031\\
2600 GA Delft, The Netherlands\\
{\tt F.M.Dekking@math.tudelft.nl}
}

\title{Generalized Beatty sequences and complementary triples}

\date{\today}

\def \proof{\noindent{\it Proof:\ \ }}

\def \endpf{{\ \ $\Box$ \medbreak}}

\newcommand{\Id}{ {\rm Id}}

\newtheorem{theorem}{Theorem}
\newtheorem{lemma}[theorem]{Lemma}
\newtheorem{corollary}[theorem]{Corollary}
\newtheorem{proposition}[theorem]{Proposition}
\theorembodyfont{\rm}

\newtheorem{remark}[theorem]{Remark}
\newtheorem{example}[theorem]{Example}
\newtheorem{definition}[theorem]{Definition}

\begin{document}

\maketitle

\begin{abstract}
A generalized Beatty sequence is a sequence $V$ defined by
$V(n)=p\lfloor{n\alpha}\rfloor+qn +r$, for $n=1,2,\dots$, where $\alpha$ is a real number,
and $p,q,r$ are integers. These occur in several problems, as for instance in  homomorphic
embeddings of Sturmian languages in the integers. Our results are for the case that
$\alpha$ is the golden mean, but we show how some results generalize to arbitrary
quadratic irrationals. We mainly consider the following question: For which sixtuples of
integers  $p,q,r,s,t,u$ are the two sequences $V=(p\lfloor{n\alpha}\rfloor+qn +r)$ and
$W=(s\lfloor{n\alpha}\rfloor+tn +u)$ complementary sequences?

We also study complementary triples, i.e., three sequences
$V_i=(p_i\lfloor{n\alpha}\rfloor+q_in+r_i), \:i=1,2,3$, with the property that the sets they
determine are disjoint with union the positive integers.

\bigskip


\bigskip


\end{abstract}

\section{Introduction}

A Beatty sequence is the sequence $A = (A(n))_{n \geq 1}$, with $A(n)=\lfloor{n\alpha}\rfloor$ for
$n \geq 1$, where $\alpha$ is a positive real number. What Beatty observed is that when
$B = (B(n))_{n \geq 1}$ is the sequence defined by $B(n)=\lfloor{n\beta}\rfloor$, with $\alpha$
and $\beta$ satisfying
\begin{equation}\label{eq:Beat}
\frac1{\alpha}+\frac1{\beta}=1,
\end{equation}
then $A$ and $B$ are \emph{complementary} sequences, that is, the sets
$\{A(n): n \geq 1\}$ and $\{B(n): n\geq 1\}$ are disjoint and their union is the set of positive
integers. In particular if $\alpha = \varphi = \frac{1 + \sqrt{5}}{2}$ is the golden ratio, this gives that
the sequences $(\lfloor n \varphi \rfloor)_{n \geq 1}$ and $(\lfloor n \varphi^2 \rfloor)_{n \geq 1}$
are complementary.

\bigskip

Among the numerous results on Beatty sequences, a paper of Carlitz, Scoville and Hoggatt
\cite[Theorem~13, p.\ 20]{car-sco-hog} studies the monoid generated by
$A=(A(n))_{n \geq 1}$ and $B=(B(n))_{n \geq 1}$ for the composition of sequences in the
case where $\alpha$ is equal to $\varphi = \frac{1 + \sqrt{5}}{2}$, the golden ratio.
(The composition of two integer sequences $U=(U(n))_{n \geq 1}$ and
$V=(V(n))_{n \geq 1}$ is the sequence $UV := U\circ V = (U(V(n)))_{n \geq 1}$, so that
the monoid generated by $A$ and $B$ is composed of sequences like $A^kB^jA^{\ell}...$,
where $A^k = AA\ldots A$ is the composition of $k$ sequences equal to $A$.)

\begin{theorem}[Carlitz-Scoville-Hoggatt]\label{thm:Car}
Let $U = (U(n))_{n \geq 1}$ be a composition of the sequences
$A = (\lfloor n \varphi \rfloor)_{n \geq 1}$ and
$B = (\lfloor n \varphi^2 \rfloor)_{n \geq 1}$, containing $i$ occurrences of
$A$ and $j$ occurrences of $B$, then for all $n \geq 1$
$$
U(n) = F_{i + 2j}  A(n) + F_{i + 2 j - 1}\, n - \lambda_U,
$$
where $F_k$ are the Fibonacci numbers ($F_0 = 0$, $F_1 = 1$, $F_{n+2} = F_{n+1} + F_n$)
and $\lambda_{ U}$ is a constant.
\end{theorem}

\medskip

\noindent {\bf Definition} We call {\em generalized Beatty sequence} any sequence $V$ of the type  $V(n) = p(\lfloor n \alpha \rfloor) + q n +r $, $n\ge 1$,  where $\alpha$ is a real number, and $p,q,$ and $r$ are integers.

\medskip

 Two examples of generalized Beatty sequences are $U=AA$ and $U=AB$, where Theorem~\ref{thm:Car} gives $AA(n)=A(n)+n-1$, and $AB(n)=2A(n)+n$. These two formulas directly imply the following  result.

\medskip

\begin{corollary}\label{corr:VA} Let $V$ be a generalized Beatty sequence given by $V(n) = p(\lfloor n \varphi \rfloor) + q n +r $, $n\ge 1$. Then $VA$ and $VB$ are generalized Beatty sequences with parameters $(p_{V\!A},q_{V\!A},r_{V\!A})=(p+q,p,r-p)$ and $(p_{W\!A},q_{W\!A},r_{W\!A})=(2p+q,p+q,r)$.
\end{corollary}

\bigskip

As an extension of Beatty's observation  the following natural questions can be asked.

\bigskip

\noindent
{\bf Question~1} \ Let $\alpha$ be an irrational number, and let $A$ defined by
$A(n)=\lfloor{n\alpha}\rfloor$ for $n\geq 1$ be the Beatty sequence of $\alpha$.
Let $\Id$ defined by $\Id(n)=n$ be the identity map on the integers.
For which sixtuples of integers $p,q,r,s,t,u$ are the two sequences \, $V = pA + q\,\Id + r$  and \, $W = s A + t\,\Id + u$ complementary sequences?

\bigskip

\noindent
{\bf Question~2} \  For which nonuples of integers  $(p_1,q_1,r_1,p_2,q_2,r_2,p_3,q_3,r_3)$  are the
three sequences $V_i = p_i A + q_i\,\Id + r_i, \:i=1,2,3$
 a complementary triple, i.e., the sets they determine are disjoint with union the positive
integers\footnote{And when is this partition ``nice''? (in Fraenkel's terminology
\cite{fraenkel1994},  a ``nice'' integer DCS --Disjoint Covering System).}.

\begin{remark}
The theorem of Carlitz, Scoville and Hoggatt above was rediscovered by Kimberling
\cite[Theorem 5, p.\ 3]{kimberling}: it is thus attributed to Kimberling in, e.g.,
\cite[p.~575]{fraenkel2010-1}, \cite[p.~647]{fraenkel2010-2}, \cite[p.~20--21]{larsson-et-al}.
This was corrected in \cite[Theorem~2, p.~2]{ballot}. Theorem~1 in \cite{griff2015} is also a
special case of the theorem of Carlitz, Scoville and Hoggatt.
\end{remark}

\begin{remark}
The four papers \cite{Mercer}, \cite{A-Fraenkel-S}, \cite{Kimberling-BG} and
\cite{Hildebrand} consider different generalizations of Beatty sequences.
\end{remark}

\begin{remark}
One can ask whether the monoid generated by other complementary sequences by
composition can be written as a subset of the set of linear combinations of a finite number of
elements. Some answers for Beatty sequences can be found in the rich paper of Fraenkel
\cite{fraenkel1994} (see, e.g., p.\ 645). Another, possibly unexpected, example is given by
the Thue-Morse sequence. Namely call odious (resp.~evil) the integers whose binary
expansion contains an odd (resp. even) number of $1$'s, then it was proved in
\cite[Corollaries~1 and 3]{all-clo-she} that the sequences
$(A(n))_{n \geq 0}$ and $(B(n))_{n \geq 0}$ of odious and evil numbers satisfy for all $n$
\begin{align*}
 &A(n) = 2n + 1 - t(n), & B(n) = 2n + t(n),\quad &A(n) - B(n) = 1 - 2 t(n) \\
 &A(A(n))  = 2A(n),     & B(B(n)) = 2B(n), \quad  &A(B(n)) = 2B(n) + 1, \quad
    B(A(n)) = 2A(n) + 1.
\end{align*}
where $(t(n))_{n \geq 0}$ is the Thue-Morse sequence, i.e., the characteristic function of
odious integers. (This sequence can be defined by $t(0) = 0$ and for all $n \geq 0$,
$t(2n) = t(n)$ and $t(2n+1) = 1 - t(n)$.) This easily implies that any finite composition of
$(A(n))_{n \geq 0}$ and $(B(n))_{n \geq 0}$ can be written as
$(\alpha A(n) + \beta B(n) + \gamma)_{n \geq 0}$, since $t(A(n)) = 1$ and $t(B(n)) = 0$
for all $n$.
\end{remark}

\section{Complementary pairs}

Let $\alpha$ be an irrational number, and let $A$ defined by $A(n)=\lfloor{n\alpha}\rfloor$
for $n\ge 1$ be the Beatty sequence of $\alpha$. In this section we consider Question~1
of the Introduction, which we call the Complementary pair problem.

\medskip

In what follows we will require that as a function $A:\mathbb{N}
\rightarrow \mathbb{N}$ is injective, since we then have a $1$-to-$1$ correspondence between
sequences and subsets of $\mathbb{N}$. (See \cite{Kim-Sto} for non-injective Beatty sequences.)

In the case that $V$ and $W$ are increasing, we will also require, without loss of generality, that
$V(1)=1$. Solutions $(p,q,r,s,t,u)$ with $p=0$ or $s=0$ will be called {\it trivial}.

\medskip

\noindent The homogeneous Sturmian sequence generated by a real number
$\alpha\in (0,1)$  is the sequence
$$c_\alpha:=(\lfloor(n+1)\alpha\rfloor-\lfloor n\alpha\rfloor)_{n\ge 1}.$$
(For more about Sturmian sequences, the reader can consult, e.g., \cite[Ch.~2]{lothaire}.)

\noindent A real number $\alpha$ is called  a {\em Sturm number} if $\alpha\in (0,1)$ is a
quadratic irrational number with algebraic conjugate $\overline{\alpha}$ satisfying
$\overline{\alpha}\notin (0,1)$. Sturm numbers have a property that is useful to recognize
their generalized Beatty sequences.

%%\section{Appendix: Lemma $4^+$ from All\&Dekk GBS}

\begin{proposition}\label{lem:3iff} {\bf (\cite{Crisp-et-al}, \cite{Allauzen})} Let $\alpha$ be a Sturm number. Then there exists a morphism $\sigma_\alpha$ on the alphabet $\{0,1\}$, such that $\sigma_\alpha(c_\alpha)=c_\alpha$.
\end{proposition}

\noindent In the following we will consider the variants of  $\sigma_\alpha$ on various other alphabets than $\{0,1\}$, but will not indicate this in the notation. The  following lemma is  implied trivially by
$$V = p A + q \, \Id + r \ \ \Rightarrow \ \ V(n+1) - V(n) = p(A(n+1) - A(n))+q = p\,c_\alpha(n)+q.$$

\begin{lemma}\label{lem:diffgen} Let $\alpha$ be a Sturm number. Let $V = (V(n))_{n \geq 1}$ be the generalized Beatty sequence defined
by $V(n) = p(\lfloor n \alpha \rfloor) + q n +r $, and let $\Delta V$ be the sequence of its first
differences. Then $\Delta V$ is the fixed point of $\sigma_\alpha$ on the alphabet $\{q, p+q\}.$
\end{lemma}

 We remark that it can be shown that the first letters of $\sigma_\alpha(0)$ and $\sigma_\alpha(1)$ are equal (see, e.g., \cite{Dekking-INT-2018}), so $\sigma_\alpha$ has a unique
 fixed point. It is also obvious that this fixed point starts with $0$ if $\alpha\in (0,1/2)$, and
 with $1$ if $\alpha\in (1/2,1)$.
 For general $\alpha$, one replaces $\alpha$ with
$\breve{\alpha}=\alpha- \lfloor \alpha \rfloor ( = \{\alpha\} )$.
When $\alpha$ is the golden mean $\varphi=(1+\sqrt{5})/2$,  the morphism generating the
sequence associated to the Sturm number $\breve{\varphi}=\varphi-1$ is $0\mapsto 1,
1\mapsto 10$, so one has to exchange $0$ and $1$ if one wishes to compare $\Delta V$
with the classical Fibonacci morphism $0\mapsto 01, 1\mapsto 0$.
 As a special case of Lemma~\ref{lem:diffgen} we therefore obtain one direction of
 the following lemma.

\begin{lemma}\label{lem:diff} Let $V = (V(n))_{n \geq 1}$ be the generalized Beatty
sequence defined by $V(n) = p(\lfloor n \varphi \rfloor) + q n +r$, and let $\Delta V$ be the
sequence of its first differences. Then $\Delta V$ is the Fibonacci word on the alphabet
$\{2p+q, p+q\}$. Conversely, if $x_{ab}$ is the Fibonacci word on the alphabet
$\{a,b\}$,  then any $V$ with $\Delta V= x_{ab}$ is a generalized Beatty sequence
$V=((a-b) \lfloor n \varphi \rfloor)+(2b-a)n+r)$ for some integer $r$.

\end{lemma}

Another observation is that the $q \, \Id+r$ part in a generalized Beatty sequence generates
arithmetic sequences. The following lemma, which will be useful in proving
Theorem~\ref{thm:sol}, shows that in some weak sense the Wythoff part $pA$ of a g
eneralized  Beatty sequence is  orthogonal to its arithmetic sequence part, provided that
$\frac13<\{\alpha\}<\frac23$, where $\{\alpha\}=\alpha-\lfloor \alpha\rfloor$. We prove this for
$\frac43<\alpha<\frac53$.

\begin{lemma}\label{lem:orth} Let $\alpha$ satisfy $\frac43<\alpha<\frac53$, and let
$V = (V(n))_{n\ge 1}$  be the generalized Beatty sequence
defined by $V(n) = p(\lfloor n \alpha\rfloor)+qn+r$ with $p\ne 0$, then neither
$(V(1),V(2),V(3))$, nor $(V(2),V(3), V(4))$
can be an arithmetic sequence of length $3$.
\end{lemma}

\proof When $\frac32<\alpha<\frac53$, we have $\lfloor\alpha\rfloor=1, \lfloor2\alpha\rfloor=3$, and $\lfloor3\alpha\rfloor=4, \lfloor4\alpha\rfloor=6$, so
\begin{eqnarray*}
V(2)-V(1) &=& p\lfloor 2\alpha \rfloor +2q+r -p\lfloor\alpha\rfloor-q-r = 2p+q,\\
V(3)-V(2) &=& p\lfloor 3\alpha \rfloor +3q+r -p\lfloor 2\alpha \rfloor-2q-r = p+q,\\
V(4)-V(3) &=& p\lfloor 4\alpha \rfloor +4q+r -p\lfloor 3\alpha \rfloor-3q-r = 2p+q,
\end{eqnarray*}
and the result follows, since $p\ne 0$. When $\frac43<\alpha<\frac32$,  we have $\lfloor\alpha\rfloor=1, \lfloor2\alpha\rfloor=2, \lfloor3\alpha\rfloor=4$, and $\lfloor4\alpha\rfloor=5$. So this time $V(2)-V(1)=p+q, V(3)-V(2)=2p+q$  and $V(4)-V(3)=p+q$, leading to the same conclusion.  $\hfill\Box$

\begin{remark}\label{comput}
We note for further use that solving the equations in the proof of Lemma~\ref{lem:orth} for $p$ and $q$, supplemented with an equation for $r$, yields in the case $\frac32<\alpha<\frac53$ that\\[-.2cm]
$$
\left\{
\begin{array}{ll}
p &= -V(1)+2V(2)-V(3)\\
q &= \;V(1) - 3V(2) + 2V(3) \\
r &= \;V(1) + V(2) - V(3).
\end{array}
\right.
$$
\end{remark}

\medskip

Let $\alpha=\varphi$  be the golden mean. Then the  classical solution is $(p,q,r)=(1,0,0)$ and
$(s,t,u)=(1,1,0)$, which corresponds to the Beatty pair
$(\lfloor n\varphi \rfloor), (\lfloor n\varphi^2 \rfloor)$. Another solution is given by
$$(p,q,r)=(-1,3,-1), \quad (s,t,u)=(1,2,0),$$
which corresponds to the Beatty pair
$(\lfloor n(5-\sqrt{5})/2 \rfloor), (\lfloor n(5+\sqrt{5})/2 \rfloor)$, which is equal to
$$(\lfloor n(3-\varphi) \rfloor), (\lfloor n(\varphi+2) \rfloor).$$

\begin{theorem}\label{thm:sol} Let $\alpha=\varphi$. Then there are exactly two non-trivial increasing
solutions to the complementary pair problem:  $(p,q,r,s,t,u)=(1,0,0,1,1,0)$ and
$(p,q,r,s,t,u)=(-1,3,-1,1,2,0)$.
\end{theorem}

\proof  Recall that $V(1)=1$. Note that  $V(2) < 5$, since otherwise $(W(1), W(2), W(3))=(2,3,4)$,
which is not allowed by Lemma~\ref{lem:orth}. There are therefore three cases to consider, according to
the value of $V(2)$.

\begin{enumerate}

\item{The case }  $V(1) = 1, V(2) = 2$. Then by Lemma~\ref{lem:orth}, $V(3) = 3$ is not
possible.\\[-.7cm]

\begin{enumerate}
\item{}
If $V(3) = 4$, then, by Remark~\ref{comput}, $p=-1$, $q=3$, $r=-1$, which is one of the two
solutions.
\item{}
If $V(3) = 5$, then, by Remark~\ref{comput}, $p=-2$, $q=5$, $r=-2$, which implies that
$V(4) = 6$, $V(5) = 7$, $V(6) = 10$. So $W(1) = 3$, $W(2) = 4$, $W(3) = 8$,
which gives  $s=-3$, $t=7$, $u=-1$ (Remark~\ref{comput} applied to $W$), implying $W(5) = 10$,
which contradicts complementarity.
\item{}
If $V(3) = m$ with $m>5$, then $W(1) = 3$, $W(2) = 4$, $W(3) = 5$, which contradicts
Lemma~\ref{lem:orth}.
\end{enumerate}

\item{The case} $V(1) = 1$, $V(2) = 3$.\\[-.7cm]

\begin{enumerate}

\item{ }
If $V(3) = 4$, then, by Remark~\ref{comput}, $p=1$, $q=0$, $r=0$, which is one of the two
solutions.

\item{ }
If $V(3) = 5$, then we obtain a contradiction with Lemma~\ref{lem:orth}.

\item{ }If $V(3) = 6$, then, by Remark~\ref{comput}, $p=-1$, $q=4$, $r=-2$, which implies
$V(5) = 10$. But we must then have $W(1) = 2, W(2) = 4, W(3) = 5$, so
(Remark~\ref{comput} applied to $W$), $s=1$, $t=0$, $u=1$, which implies $W(6) = 10$,
a contradiction with complementarity.

\item{ }
If $V(3)=m$ with $m>6$, then we obtain a contradiction with Lemma~\ref{lem:orth},
since then $W(2) = 4$, $W(3) = 5$, $W(4) = 6$.

\end{enumerate}

\item{The case} $V(1) = 1$, $V(2) = 4$.\\[-.7cm]

\begin{enumerate}

\item{ }
If $V(3) = 5$, then, by Remark~\ref{comput}, $p=2$, $q=-1$, $r=0$, thus $V(4) = 8$;
hence $W(1) = 2$, $W(2) = 3$, $W(3) = 6$. Hence, by Remark~\ref{comput} applied to $W$, $s=-2$,
$t=5$, $u=-1$, so that $W(5) = 8 = V(4)$, which contradicts complementarity.

\item{ }
If  $V(3) = 6$, then $W(1) = 2$, $W(2) = 3$, $W(3) = 5$.
Thus, by Remark~\ref{comput} applied to $W$, $s=-1$, $t=3$, $u=0$. Hence $W(4) = 6 = V(3)$,
which contradicts complementarity.

\item{ }
If $V(3) = 7$, then we obtain a contradiction with Lemma~\ref{lem:orth}.

\item{ }
If $V(3) = m$ with $m>7$, then it follows that  $V(3) = 8$, since  we have $W(1) = 2$,
$W(2) = 3$, $W(3) = 5$, yielding, by Remark~\ref{comput} applied to $W$,
$W(n) = (-A(n) + 3n) = 2, 3, 5, 6, 7, 9, 10, 12, 13, 14, \dots$ With $V(3) = 8$, one obtains
(by Remark~\ref{comput}) that $V(n)=-A(n)+5n-3$, but then $V(5) = 14 = W(10)$, i.e., $V$
and $W$ are not complementary. \endpf
\end{enumerate}
\end{enumerate}

For $\alpha=\sqrt{2}$\, the classical solution to the complementary pair problem is $V=A$,
$W=A+2\,\Id$, i.e., the Beatty pair given by $V(n)=\lfloor n\sqrt{2}\rfloor$, and
$W(n)=\lfloor n(2+\sqrt{2})\rfloor$. As $\frac43<\sqrt{2}<\frac32$, we can use
Lemma~\ref{lem:orth} and adapt Remark~\ref{comput} to prove the following result, in the
same way as Theorem~\ref{thm:sol}.

\begin{theorem}\label{thm:pell1} Let $\alpha=\sqrt{2}$. Then there is a unique non-trivial
increasing solution to the complementary pair problem:  $(p,q,r,s,t,u)=(1,0,0,1,2,0)$.
\end{theorem}

We end this section with an example where $\{\alpha\}\notin (\frac13,\frac23)$.

\begin{theorem}\label{thm:pell2} Let $\alpha=\sqrt{8}$. Then there is a unique non-trivial increasing
solution to the complementary pair problem:  $(p,q,r,s,t,u)=(1,4,0,-1,4,0)$.
\end{theorem}

\proof Since $(4+\sqrt{8}, 4-\sqrt{8})$ is a Beatty pair, $(p,q,r,s,t,u)=(1,4,0,-1,4,0)$ is a solution to the complementary pair problem. To prove that it is unique is more involved. We fix $V(1)=1$.

Let $\breve{\alpha}=\alpha-2=\sqrt{8}-2$. Then $\breve{\alpha}\in (0,1)$, and $\breve{\alpha}$ has the periodic continued fraction expansion $[0;\overline{1,4}]$. It follows then from \cite{Crisp-et-al}, or from Corollary 9.1.6 in \cite{AllShall} that the morphism $\sigma_{\breve{\alpha}}$ fixing the homogeneous Sturmian sequence $c_{\breve{\alpha}}$ is given by
$$   \sigma_{\breve{\alpha}}: 0 \mapsto11110, \quad    1 \mapsto 111101.$$
Note that
$$V(n)=p\lfloor n\sqrt{8}\rfloor+qn+r= p\lfloor n(\sqrt{8}-2)\rfloor+(2p+q)n+r=
p\lfloor n\breve{\alpha}\rfloor+(2p+q)n+r.$$
The difference sequence $\Delta V$ of $V$ is therefore the fixed point of  $\sigma_{\breve{\alpha}}$ on the alphabet $\{2p+q, 3p+q\}$. Since we require $V$ to be increasing, both $2p+q$ and $3p+q$ have to be larger than 0.
 We split the  possibilities according to the value of $3p+q$. The arguments below are based on the fact, following from the form of $\sigma_{\breve{\alpha}}$, that $V$ starts with an arithmetic sequence of length 5, followed by an arithmetic sequence of length 6, both with common differences $3p+q$, and separated by a distance $2p+q$.

\begin{enumerate}

\item{The case } $3p+q\ge 3$.\\[-.5cm]

If  $3p+q\ge 3$, then $W(1)=2,\, W(2)=3$, so $W=sA+t\Id+u$ has to start with an arithmetic sequence of length 5 with common difference $1$, i.e., $W(1),W(2),\dots,W(5) = 2,3,\dots,6$. Moreover, since $W(6)=7$ is not possible (it would imply $p=0$), $V(2)=7$, which implies $V(3)=13$. But then the second arithmetic sequence of $W$, which has length 6, does not fit in between $V(2)$ and $V(3)$.

\item{The case } $3p+q=2$.\\[-.5cm]

In this case $V(1),\dots ,V(5) = 1,3,5,7,9$, so $W(1),\dots ,W(5) =2,4,6,8,10$. Then either $V(6)=11$, or $W(6)=11$.

In the former case we must have $2p+q=V(6)-V(5)=2$, which implies  $p=0$, which is trivial.

In the latter case $W(6),\dots,W(11)=11,13,15,17,19,21$, and $W(12)=22$, since $W(12)-W(11)=W(6)-W(5)=1$. But also, $V(6),\dots,V(11)$ equals $12,14,\dots,22$. So $V$ and $W$ are not complementary.

\item{The case } $3p+q=1$.\\[-.5cm]

In this case $V(1),\dots ,V(5) = 1,2,3,4,5$, so $W(1)=6$, since $V(6)=6$ would imply $p=0$.
Then either $V(6)=7$, or $W(2)=7$.

In the former case, $V(6),\dots ,V(11) = 7,8,9,10,11,12$, and $W(2)=13$. This implies that
$2=V(6)-V(5)=2p+q$, which leads to $(p,q,r)=(-1,4,0)$, and $(s,t,u)=(1,4,0)$, which is the
announced solution.

In the latter case $W(1),\dots ,W(5) =6,7,8,9,10$, and $V(6)=11$. So $2p+q=V(6)-V(5)=5$.
This implies $(p,q,r)=(-5,16,-5)$, and $(s,t,u)=(-6,19,-1)$. But then $V(12)=22$, and
$W(11)=22$. So $V$ and $W$ are not complementary.  \endpf
\end{enumerate}

\subsection{Generalized Pell equations}
If $V$ and $W$ are not increasing, then an analysis as in the proof of
Theorem~\ref{thm:sol} is still possible, but very lengthy. We therefore consider another
approach in this subsection. Considering the densities of $V$ and $W$ in $\mathbb{N}$,
one sees that a \emph{necessary} condition for $(pA+q\,\Id+r,  sA+t\,\Id+u)$ to be a
complementary pair is that
\begin{equation}\label{eq:dens}
\frac1{p\alpha+q}+\frac1{s\alpha +t }=1
\end{equation}

In what follows we concentrate on the case $\alpha = \varphi = (1+\sqrt{5})/2$, but our
arguments can be generalized to the case of arbitrary quadratic irrationals.

\smallskip

\begin{proposition}\label{necessary-cpp}
A necessary condition  for the pair $V=pA+q\,\Id+r \;\;{\rm and}\;\; W=s A+t\,\Id+u$ to be a
complementary pair is that $p\ne 0$ is a solution to the generalized Pell equation
$$5p^2 x^2- 4x =y^2, \quad x,y\in \mathbb{Z}$$.
\end{proposition}
\vspace*{-.7cm}
\proof
Using $\varphi^2=1+\varphi$, a straightforward manipulation shows that \eqref{eq:dens} implies
\begin{equation*}
(ps+pt+qs-p-s)\varphi=q+t-ps-qt.
\end{equation*}
But since $\varphi$ is irrational, this can only hold if
\begin{equation}\label{eq:cond}
ps +pt+qs-p-s=0, \quad q+t-ps-qt=0.
\end{equation}
The first equation gives $pt=p-(p+q-1)s$. Eliminating $pt$ from $p^2s+(q-1)pt-pq=0$, we obtain
$p^2s+(p-(p+q-1)s)(q-1)-pq=0.$
This gives the quadratic equation $$s\,q^2+(p-2)s\,q-(p^2+p-1)s+p=0.$$
Since $q$ is an integer, $\Delta:=(p-2)^2s^2+4s((p^2+p-1)s-p)$ has to be an integer squared. Trivial manipulations yield that
\begin{equation}\label{eq:Delta}\Delta=5p^2s^2-4ps.\end{equation}
 Since $p$ divides the square $\Delta$, $5p^2s^2-4ps=p^2y^2$ for some integer $y$, and hence $p$ also divides $s$.
If we put $s=px$, we obtain  $5p^3x^2-4p^2x=p^2y^2$, which finishes  the proof of the proposition.
$\hfill\Box$

\bigskip

Actually there is a simple characterization of the integers $p$ such that the Diophantine equation
above has a solution.

\begin{proposition}\label{general-pell}
The generalized Pell equation
$$
5p^2 x^2- 4x =y^2, \quad x,y\in \mathbb{Z}
$$
has a solution for $p > 0$ if and only if $p$ divides some Fibonacci number of odd index, i.e., if and
only $p$ divides some number in the set $\{1, 2, 5, 13, 34, \ldots\}$.
\end{proposition}

\proof
First suppose that there are integers $p > 0$ and $x,y \in  \mathbb{Z}$ such that
$5p^2 x^2 - 4x =y^2$. Let $d := \gcd(x,y)$ and $x' = x/d$, $y'=y/d$, so that $\gcd(x',y') = 1$.
We thus have
$$
5 p^2 d x'^2 - 4 x' = d y'^2.
$$
Thus $x'$ divides $dy'^2$, but it is prime to $y'$, hence $x'$ divides $d$. Since clearly $d$ divides
$4x'$, we have $d = \alpha x'$ for some $\alpha$ dividing $4$, hence $\alpha$ belongs to
$\{1, 2, 4\}$. This yields $\alpha(5p^2 x'^2 - y'^2) = 4$. We distinguish three cases.
\,
\begin{enumerate}
\,
\item{  } If $\alpha = 1$, then we have $5p^2 x'^2 - y'^2 = 4$. But the equation $5X^2 - 4 = Y^2$
has an integer solution if and only if $X$ is a Fibonacci number with odd index \cite[p.~91]{Lind}.
Hence $px'$ must be a Fibonacci number with odd index, thus $p$ divides a Fibonacci number
with odd index.
\,
\item{  } If $\alpha = 2$, then we have $5p^2 x'^2 - y'^2 = 2$. Note that $x'$ must be odd, otherwise
$x'$ and $y'$ would be even, which contradicts $\gcd(x',y') = 1$. Thus $5p^2x'^2 \equiv p^2 \bmod 4$,
hence $p^2 - 2 \equiv y'^2 \bmod 4$. If $p$ is even, this yields $y'^2 \equiv 2 \bmod 4$, while if
$p$ is odd, this gives $y'^2 \equiv 3 \bmod 4$. There is no such $y'$ in both cases.
\,
\item{  } If $\alpha = 4$, then we have $5p^2 x'^2 - y'^2 = 1$, thus $5(2px')^2 - (2y')^2 = 4$, then
$2px'$ must be a Fibonacci number with odd index, thus $p$ divides a Fibonacci number with odd
index.
\,
\end{enumerate}
\,
Now suppose that $p$ divides some Fibonacci number with odd index, say there exists a $k$ with
$F_{2k+1} = p\beta$. We will construct an integer solutions in $(x,y)$ to the equation
$5p^2 x^2 - 4x = y^2$. We know (again \cite[p.~91]{Lind}) that there exists some integer $\gamma$
with $5 F_{2k+1}^2 - 4 = \gamma^2$ thus $5p^2 \beta^2 - 4 = \gamma^2$. Let $x = \beta^2$ and
$y = \beta \gamma$. Then
$$
5 p^2 x^2 - 4 x = 5 p^2 \beta^4 - 4 \beta^2 = \beta^2(5p^2 \beta^2 - 4) = \beta^2 \gamma^2 = y^2.
\ \ \hfill\Box
$$
\begin{corollary} \label{legendre}
There are no solutions to the complementary pair problem if $-1$ is not a square modulo $p$, i.e.,
if $p$ does not belong to the sequence $1, 2, 5, 10, 13, 17, 25, 26, 29, 34, 37, 41, \ldots$ (sequence
A008784 in \cite{oeis}). This is in particular the case if $p$ has a prime divisor congruent to $3$
modulo $4$.
\end{corollary}
\,
\proof
We will prove that if there are solutions to the complementary problem for $p$,
thus if $p$ divides an odd-indexed Fibonacci number (Propositions~\ref{necessary-cpp} and
\ref{general-pell}), then $-1$ is a square modulo $p$. Using again the characterization in
\cite[p.~91]{Lind}, there exist two integers $x,y$ with $5p^2x^2 - 4 = y^2$. We distinguish two cases.
\begin{enumerate}
\,
\item{  } If $p$ is odd, we have $y^2 \equiv -4 \bmod p$ and $2^2 \equiv 4 \bmod p$. But $2$
is invertible modulo $p$, hence, by taking the quotient of the two relations, we obtain that $-1$
is a square modulo $p$.
\,
\item{  } If $p$ is even, remembering that $px = F_{2k+1}$ for some $k$, we claim that $p$ must be
congruent to $2$ modulo $4$ and that $x$ must be odd. Namely the sequence of odd-indexed
Fibonacci numbers, reduced modulo $4$, is easily seen to be the periodic sequence
$(1 \ 2 \ 1)^{\infty}$. Hence it never takes the value $0$ modulo $4$. The equality $5p^2x^2 - 4 = y^2$
implies that $y$ must be even, thus we have $5(p/2)^2 x^2 - 1 = (y/2)^2$, say $(y/2)^2 = -1 + z(p/2)$.
Up to replacing $(y/2)$ with $(y+p)/2$, we may suppose that $(y/2)$ is even (recall that
$p/2$ is odd). Thus $z(p/2)$ is even, hence $z$ is even, say $z = 2z'$. This gives $(y/2)^2 = -1 + z'p$,
thus $-1$ is a square modulo $p$.
\,
\end{enumerate}
\,
\begin{remark}
We have just seen that if the integer $p$ divides some odd-indexed Fibonacci number
then $-1$ is a square modulo $p$ (sequence A008784 in \cite{oeis}). A natural question is then
whether it is true that if $-1$ is a square modulo $p$, then $p$ must divide some odd-indexed
Fibonacci number. The answer is negative, since on one hand $12^2 \equiv - 1 \bmod 29$, and, on
the other hand, the sequence of odd-indexed Fibonacci numbers modulo 29 is the periodic sequence
$(1, 2, 5, 13, 5, 2, 1)^{\infty}$ which is never zero.
\end{remark}

\medskip

Let us look at examples of solutions to the Diophantine equation for values of $p$ that divide
some Fibonacci number with odd index.
Consider, for example, the case where $p=s$. Then Equation \eqref{eq:Delta} becomes
$\Delta=5p^4-4p^2$, so the Diophantine equation is
$$5x^2-4=y^2, \quad x,y \in \mathbb{Z}.$$
%It is well known that it has as solutions the Fibonacci numbers with odd indices, i.e., we have found
%that necessarily $p$ has to take a value  in the set $\{1,2,5,13,34,\dots\}$.

For $p=F_1=1$ we obtain the two sequences $V=A+r$ and $W=A+\Id+u$. These are
complementary only when $r=u=0$, and we obtain the classical Beatty pair $(A,A+\Id)$.

For $p=F_3=2$ we obtain the two sequences $V=2A+2\,\Id+r$ and $W=2A-2\,\Id+u$. These
cannot be complementary for any $r$ and $u$, since for $u=0$ we have
$W(n)= 2(\lfloor n \varphi\rfloor) - 2n= 2(\lfloor n(\varphi-1)\rfloor$, which gives all even numbers,
since $\varphi-1<1$. This an example where Equation \eqref{eq:dens} does not apply, since $W$
as a function is not injective.

For $p=F_5=5$ we obtain the two sequences $V=5A+4\,\Id+r$ and $W=5A-7\,\Id+u$.
To make these complementary we are forced to choose $r=u=3$,  and we obtain

$V=(12, 26, 35, 49, 63, 72, 86, 95, 109, 123, 132, 146, 160, 169, 183, 192, 206, 220, 229, 243, 252, 266, \dots),$

$W=(1, 4, 2, 5, 8, 6, 9, 7, 10, 13, 11, 14, 17, 15, 18, 16, 19, 22, 20, 23, 21, 24, 27, 25, 28, 31, 29, 32, 30, \dots).$

\noindent Now a proof that $V$ and $W$ form a complementary pair is much harder, when
we let $V$ start with $V(0)=3$, to include $3$ in the union. We can perform the following
trick. We split $W$ into $(W(A(n)))_{n \geq 1}$, and $(W(B(n)))_{n \geq 1}$
(cf. Corollary~\ref{corr:VA}). The two sequences $WA$ and $WB$ are increasing, and we
can prove that $(V(n))_{n \geq 0}$, $(W(A(n)))_{n \geq 1}$, and $(W(B(n)))_{n \geq 1}$ form
a partition of the positive integers by proving that the  three-letter sequence obtained by
applying the morphism $0 \mapsto 1120,  1\mapsto 11100$ to  the fixed point of the
morphism $g$ given by  $g: 0 \mapsto 01, \ 1 \mapsto 011$,  has the property  that the
preimages of $0$, $1$ and $2$ are precisely these three sequences.
See Theorem~\ref{sequence-v} and its proof for a similar result.

 For $p=F_{2m+1}\ge 13$ it seems that we can always choose $r$ and $u$ for  in such a
 way that we get almost complementary sequences:
namely, e.g., for $p=13$ we find $q=9$ and $t=-20$. If we take $r=u=9$, then we
\underline{almost} get a complementary pair. One finds
$V = 9,31,66,88,123,158,180,215,\dots$ and $W=2, 8, 1, 7, 13, 6, 12, 5, 11, 17, 10\dots$.
So $3$ and $4$ are missing. We thought we could prove, perhaps using something like
the Lambek-Moser Theorem \cite{Moser}, that for all $F_{2m+1} > 5$ the two sequences
are complementary excluding finitely many values, but we were not successful.

\section{Complementary triples}

Here we will find several complementary triples consisting of sequences
$V_i = p_iA+q_i\,\Id+r_i,\;\; i=1,2,3$, where $A(n)=\lfloor n\alpha \rfloor$, and $\alpha$
is a real number.

It is interesting that the case $p_1=p_2=p_3=1$ cannot be realized. This was proved by
Uspensky in 1927, see \cite{fraenkel1977}.

The case with different $\alpha_i$, is analysed by Tijdeman \cite{Tijdeman} for {\it rational}
$\alpha_i, \:1=1,2,3$. Also see \cite{Tijdeman2} for the inhomogeneous Beatty
case $(V_i(n))_n = (\lfloor n\alpha_i+\beta_i \rfloor)_n, \ i=1,2,3$.

\medskip

\noindent There is one triple in which we will be particularly interested
(see Theorem~\ref{sequence-v}):
$$((p_1,q_1,r_1),(p_2,q_2,r_2),(p_3,q_3,r_3))=((2,-1,0), (4,3,2), (2,-1,2)).$$

\medskip

\noindent We  allow that the sequences $(V_i)$ are each indexed either by $\{0,1,2,\dots\}$
or by $\{1,2,\dots\}$.

%Then  $(1,2,0),(1,2,1),(2,-1,1)$ is a solution (obtained by splitting one sequence in the complementary pair $(A,A+\Id)$.

\subsection{Two classical triples}

In this subsection $\alpha$ is always the golden mean $\varphi$. Let once more
$A(n) = \lfloor n \varphi \rfloor$ for $n \geq 1$ be the terms of the lower Wythoff sequence,
and let $B$ be given by $B(n)=\lfloor{n\varphi^2}\rfloor$ for $n\ge 1$, the upper
Wythoff sequence. Then we have the disjoint union
\begin{equation}\label{eq:pair} A(\mathbb{N}) \cup B(\mathbb{N}) =\mathbb{N}.
\end{equation}
Since $B=A+\Id$, this is the classical complementary pair ($(1,0,0),(1,1,0))$.

\medskip

Here is a way to create complementary triples from complementary pairs.

\begin{proposition}\label{prop:PtoT} Let $(V,W)$ be a golden mean complementary pair $V = p A + q\,\Id + r$ and $ W = s A + t\,\Id + u$. Then
$(V_1,V_2,V_3)$ is a complementary triple, where the three parameters of $V_1$ are
$(p+q, p, r-p)$, those of $V_2$ are $(2p+q,p+q,r)$, and $V_3=W$.
\end{proposition}

\proof
 Substituting Equation~\eqref{eq:pair} in $V(\mathbb{N}) \cup W(\mathbb{N})=\mathbb{N}$ we obtain the disjoint union
\begin{equation*}\label{eq:trip}
 V(A(\mathbb{N})) \cup V(B(\mathbb{N}))\, \cup \, W(\mathbb{N})=\mathbb{N}.
\end{equation*}
Then Corollary~\ref{corr:VA} implies the statement of the proposition.
\endpf

\begin{remark}
Actually Proposition~\ref{prop:PtoT} and Corollary~\ref{corr:VA} can be generalized to cover
an infinite family of quadratic irrationals, but  their statements will not be true for all quadratic
irrationals.    We hope to revisit this point in a future article.
\end{remark}

Applying Proposition \ref{prop:PtoT} to the basic complementary pair ($(1,0,0),(1,1,0))$ gives that
$$
((1,1,-1),(2,1,0), (1,1,0))\  \mbox{\rm and } ( (1,0,0), (2,1,-1), (3,2,0))
$$ %\\
are complementary triples\footnote{In \cite{oeis} these are (A003623, A003622, A001950)
and (A000201, A035336, A101864).}, which we will call classical triples. The first classical
triple is given at the end of  Skolem's 1957 paper \cite{Skolem}.

Let $w = 1231212312312\dots $ be the fixed point of the morphism
$$1 \mapsto 12,\: 2 \mapsto 3,\: 3 \mapsto 12.$$
Then $w^{-1}(1)=AA, \,w^{-1}(2)=B$ and $w^{-1}(3)=AB$ give the three sequences
$V_1, V_3$ and $V_2$ of the first classical triple (see \cite{Dekking-JIS}).

The question arises: is there  also a morphism generating the second triple? The answer is
positive.

\begin{proposition} Let $(V_1,V_2,V_3) = (A,\, 2A+\Id\!-\!1, \,3A+2\Id)=(A,\,BA,\,BB)$. Then
$(V_1,V_2,V_3)$ is a complementary triple. Let $\mu$  be the morphism on $\{1,2,3\}$
given by
$$\mu: 1 \mapsto 121,\: 2 \mapsto 13,\: 3 \mapsto 13,$$
with fixed point $z$. Then $z^{-1}(1) = V_1,\, z^{-1}(2) = V_2$ and $z^{-1}(3) = V_3$.
\end{proposition}

\proof
The four words of length $3$ occurring in the infinite Fibonacci word $x_{\rm \scriptstyle F}$
are $010, 100, 001, 101$.
Coding these with the alphabet $\{1,2,3,4\}$ in the given order, they generate the $3$-block
morphism $\hat{f}_3$ that describes the successive occurrences of the words of length $3$
in $x_{\rm \scriptstyle F}$ (cf.~\cite{Dekking-JIS}). It is given by
$$\hat{f}_3(1) = 12,\quad \hat{f}_3(2) = 3,\quad \hat{f}_3(3) = 14,\quad \hat{f}_3(4) = 3.$$
It has just one fixed point, which is
$$z':= 1,2,3,1,4,1,2,3,1,2,3,1,4,1,2,3,\ldots .$$
We claim that
$${z'}^{-1}(1)=AA,\, {z'}^{-1}(2)=BA, \quad {z'}^{-1}(3)=AB, \quad {z'}^{-1}(4)=BB.$$
To see this, note that the $3$-block $010$ in $\xF$ uniquely decomposes as $010=f(0)0$.
It follows that the $m^{\rm th}$ occurrence of $010$ in $\xF$ corresponds exactly to the
$m^{\rm th}$ occurrence of $0$ in $f^{-1}(\xF)=\xF$. This implies that the positions of the
occurrences of $010$ are of the form $A(A(n))$, and also that the occurrences of $001$
are  of the form $A(B(n))$, since $B(\N)$ is the complement of $A(\N)$.

For the $3$-block $100$, we note that it always occurs in $\xF$ as factor of $0100$, which
uniquely decomposes in $\xF$ as $0100=f^2(0)0$. It follows that  the $m^{\rm th}$
occurrence of $100$ in $\xF$ corresponds exactly to the  $m^{\rm th}$ occurrence of $0$ in
$f^{-2}(\xF)=\xF$. This implies that the positions of the occurrences of $100$ are of the form
$B(A(n))$, and also that the occurrences of $101$ are of the form $B(B(n))$.

Since the $0$'s in $\xF$ occur either as prefix of $001$, or of $010$, we see that we have
to merge the letters $1$ and $3$ to obtain the sequence $A$. This is not possible with
$\hat{f}_3$. However, the square of this $3$-block morphism is given by
$$1 \mapsto 123,\: 2 \mapsto 14,\: 3 \mapsto 123,\: 4 \mapsto 14,$$
and now we \emph{can} consistently merge $1$ and $3$ to the single letter $1$, obtaining
the morphism $\mu$, after mapping $4$ to $3$. Under this projection the sequence $z'$ maps to $z$.
\hfill$\Box$

\subsection{Non-classical triples}

Let $\Lan$ be a language, i.e., a sub-semigroup of the free semigroup generated by a finite
alphabet under the concatenation operation. A homomorphism of $\Lan$ into the natural
numbers is a map $\Sh:\Lan\rightarrow {\mathbb N}$ satisfying $\Sh(vw)=\Sh(v)+\Sh(w)$,
for all  $v,w \in \Lan.$

Let $\Lan_{\rm \scriptstyle F}$ be the Fibonacci language, i.e., the set of all words occurring
in the Fibonacci word $x_{\rm \scriptstyle F}$,  the iterative fixed point of the morphism
$f$ defined on $\{0, 1\}^*$ by $f: 0 \mapsto 01$, $1 \mapsto 0$.
The following result is proved in \cite{Dekking-TCS-2018}.

 \begin{theorem}\label{th:Fib}{\bf \rm (\cite{Dekking-TCS-2018})}
Let  $\Sh: \Lan_{\rm \scriptstyle F}\rightarrow \mathbb{N}$ be a homomorphism.
Define $a=\Sh(0), b=\Sh(1)$. Then  $\Sh(\Lan_{\rm \scriptstyle F}$) is the union of the two
generalized Beatty sequences $ \big((a-b)\lfloor n\varphi \rfloor+(2b-a)n\big)$ and
$ \big((a-b)\lfloor n\varphi \rfloor+(2b-a)n+a-b\big)$.
\end{theorem}

For a few choices of $a$ and $b$, the two sequences in  $\Sh(\Lan_{\rm \scriptstyle F})$
and the sequence ${\mathbb N}\setminus \Sh(\Lan)$ form a complementary triple of
generalized Beatty sequences. The goal of this section is to prove this for $a=3,\,b=1$. It
turns out that the three sequences
$$(2\lfloor n\varphi \rfloor-n)_{n\ge 1},\; (2\lfloor n\varphi \rfloor-n+2)_{n\ge 1},\;
(4\lfloor n\varphi \rfloor+3n+2)_{n\ge 0},$$
form a complementary triple. % \footnote{lien avec A276886 ?}

\begin{remark}
Note that the indices for $(4 \lfloor n \varphi \rfloor) + 3n + 2)_{n \geq 0}$ are $(n\geq 0)$,
not $(n\geq 1)$
\end{remark}

\bigskip

It is easy to see that the Fibonacci word $x_{\rm\scriptstyle F}$ can be obtained as an
infinite concatenation of two kinds of blocks, namely $01$ and $001$ (part (i) of
Lemma~\ref{kimb-seq} below). Kimberling introduced in the OEIS \cite{oeis} the sequence
A284749 obtained by replacing in this concatenation every block $001$ by $2$. We let
$x_{\rm\scriptstyle K} = A284749$ denote this sequence.

\begin{lemma}\label{kimb-seq}
Let $f$, $g$, $h$, $k$ be the morphisms defined on $\{0, 1\}^*$ by
$$
f: 0 \mapsto 01, \ 1 \mapsto 0; \ \ \
g: 0 \mapsto 01, \ 1 \mapsto 011;  \ \ \  h: 0 \mapsto 01, \ 1 \mapsto 001; \ \ \
k: 0 \mapsto 01, \ 1 \mapsto 2.
$$
Then \qquad {\rm (i)\:} $x_{\rm\scriptstyle F} = f^{\infty}(0) = h(g^{\infty}(0))$,
\quad {\rm (ii)\:} $x_{\rm\scriptstyle K}  = k(g^{\infty}(0))$.
\end{lemma}

\proof

(i) An easy induction proves that for all $n \geq 0$ one has $h g^k = f^{2k} h$. (Note that it
suffices to prove that the values of both sides are equal when applied to $0$ and to $1$.)
By letting $n$ tend to infinity this implies $h g^{\infty}(0) = f^{\infty}(0)$.

(ii) Assertion (i) gives that $x_{\rm\scriptstyle F}$ is an infinite concatenation of
blocks $h(0)=01$ and $h(1)=001$, obtained as image under $h$ from $g^{\infty}(0)$. So substituting 2 for 001 in $x_{\rm\scriptstyle F}$ is the same as substituting 01 for 0 and 2 for 1 in $g^{\infty}(0)$. \endpf

\noindent It is interesting that $x_{\rm\scriptstyle K}$   is fixed point of a morphism $i$, given by $i : 0 \mapsto 01, \ 1 \mapsto 2, \ 2 \mapsto 0122.$. This follows from the relation $k g^n = i^{n+1}$ for all $n$, which is easily proved by induction.

\begin{lemma}\label{2after0}
Define the morphism $\ell$ from $\{0, 1\}^*$ to $\{0, 1, 2\}^*$ by $\ell: 0 \mapsto 012$,
$1 \mapsto 0022$. Then the sequence $v = (v_n)_{n \geq 1} = \ell (g^{\infty}(0))$ is
obtained from $x_{\rm\scriptstyle K} $ by replacing $1$ by $0$ in all blocks
$0122$ (but not in $0120$).
\end{lemma}

\proof Note that $kg:  0 \mapsto 012, \ 1 \mapsto 0122$. Lemma \ref{2after0} then follows from  $x_{\rm\scriptstyle K}  = k(g^{\infty}(0))=kg(g^{\infty}(0))$.\endpf

%\begin{lemma}\label{seq-w}
%Let $w$ be the sequence obtained from $v$ by replacing all $2$'s
%by $1$'s. Let $m$ be the morphism defined on $\{0, 1\}^*$ by $m: 0 \mapsto 011$,
%$1 \mapsto 0011$. Then $w = m(g^{\infty}(0))$.
%\end{lemma}

%\proof Letting $q$ the morphism defined by $q: 0 \mapsto 0$, $1 \mapsto 1$,
%$2 \mapsto 1$, one has $w = q(v) = q(\ell(g^{\infty}(0))) = m(g^{\infty}(0))$ since, clearly,
%$q \ell = m$. \endpf

\begin{theorem}\label{sequence-v}
Let $v$ be the sequence defined above, i.e., $v = \ell(g^{\infty}(0))$, where $g$ and
$\ell$ are the morphisms defined by $g: 0 \mapsto 01$, $1 \mapsto 011$ and
$\ell: 0 \mapsto 012$, $1\mapsto 0022$. Then the increasing sequences of
integers defined by $v^{-1}(0)$, $v^{-1}(1)$, $v^{-1}(2)$ form a
partition of the set of positive integers ${\mathbb N}^*$. Furthermore
\begin{itemize}
\item{ }
$v^{-1}(0) = \{1, 4, 5, 8, 11, 12, 15, 16, 19, 22, \ldots\}$ is equal to the sequence of
integers $(2 \lfloor n \varphi \rfloor - n)_{n \geq 1}$, where $\varphi$ is the golden ratio
$\frac{1+\sqrt{5}}{2}$ (sequence{ \rm A050140} in {\rm \cite{oeis}}),
\item{ }
$v^{-1}(1) = \{2, 9, 20, 27, \ldots\}$ is equal to the sequence of integers
$(4 \lfloor n \varphi \rfloor + 3n + 2)_{n \geq 0}$.
\item{ }
$v^{-1}(2) = \{3, 6, 7, 10, 13, 14, 17, 18, 21, 24, \ldots\}$ is equal to the sequence of
integers $((2 \lfloor n \varphi \rfloor - n+2)_{n \geq 1})$ (i.e., $2+${\rm A050140}).
\end{itemize}
\end{theorem}

\proof We see from Lemma \ref{2after0} that the positions of 2 in $v$ are the same as the positions of 2 in $x_{\rm\scriptstyle K}$. In Section  \ref{sec:return}   it is proved that $x_{\rm\scriptstyle K}^{-1}(2) = ((2 \lfloor n \varphi \rfloor - n+2)$, see Example \ref{mapsto2}, and so the third assertion of the theorem follows.

Inspection of the occurrences of 0 and 2 in $\ell(0)$ and $\ell(1)$ then shows that the first assertion will also be true.

\medskip

\noindent For the proof of the second assertion consider $\ell(01) = 0{\bf 1}20022$,\,
$\ell(011) =  0{\bf 1}200220022.$
Since $g^{\infty}(0)$ is a concatenation of the words $g(0)=01$ and $g(1)=011$, we see
 from this that the differences of indices of the positions where $1$'s in $v$ occur are
 $7$ or $11$, and moreover, that a $0$ in $g^{\infty}(0)$ generates a difference $7$, a $1$
 in $g^{\infty}(0)$ generates a difference $11$.

It is well known and easy to prove that $g^{\infty}(0)$ equals the binary complement of the
Fibonacci word prefixed with the letter $1$. From this it follows that $\Delta v^{-1}(1)$ is the
Fibonacci word on the alphabet $\{11,7\}$. Now Lemma~\ref{lem:diff}  gives the generalized
Beatty sequence $V=(4 \lfloor n \varphi \rfloor + 3n + 2)$. The first element $2$ in
$v^{-1}(1)$ is obtained by letting $V$ start at $n=0$ instead of $n=1$.
\endpf

\begin{remark}
Some of the sequences above are images of Sturmian sequences by a morphism. Namely
$v = \ell(g^{\infty}(0))$, $x_{\rm\scriptstyle K}  = k(g^{\infty}(0))$.
Such sequences are examples of sequences called \emph{quasi-Sturmian} in \cite{Cassaigne}.
Their block complexity is of the form $n+C$ for $n$ large enough ($C=1$ for Sturmian sequences).
These sequences were studied, e.g., in \cite{Paul}, \cite{Coven}, and \cite{Cassaigne}.
\end{remark}

\section{Generalized Beatty sequences and return words}\label{sec:return}

In this  section we show that generalized Beatty sequences are closely related to return
words.

\begin{theorem}\label{th:GBS-returns}
Let $x_{\rm \scriptstyle F}$ be the Fibonacci word, and let $w$ be any  word in the
Fibonacci language $\Lan_{\rm \scriptstyle{F}}$.
Let $Y$ be the sequence of positions of the occurrences of $w$ in $x_{\rm \scriptstyle F}$.
Then $Y$ is a generalized Beatty sequence, i.e., for all $n\ge 0$,
$Y(n+1) = p\lfloor n \varphi \rfloor + qn +r$ with parameters $p,q,r$, which can be explicitly
computed.
\end{theorem}

\proof Let $x_{\rm \scriptstyle F}=r_0(w)r_1(w)r_2(w)r_3(w)\dots $, written as a
concatenation of return words of the word $w$ (cf. \cite{HuangWen}, {Lemma 1.2}).
According to Theorem2.11 in \cite{HuangWen}, if we skip $r_0(w)$, then the return words
occur as the Fibonacci word on the alphabet $\{r_1(w),r_2(w)\}$. Thus the distances
between occurrences of $w$ in $x_{\rm \scriptstyle F}$ are equal to
$l_1:=|r_1(w)|$ and $l_2:=|r_2(w)|$. We can apply Lemma~\ref{lem:diff}, which yields
$p=l_1-l_2, \; q=2l_2-l_1$. Inserting $n=0$,  we find that $r=|r_0(w)|+1$, as the first
occurrence of $w$ is at the beginning of $r_1(w)$. \endpf

\subsection{The Kimberling transform}

Here we will obtain non-classical triples appearing  in another way, namely as the
three indicator functions $x^{-1}(0),x^{-1}(1)$ and $x^{-1}(2)$, of a sequence $x$ on an
alphabet $\{0,1,2\}$ of three symbols. In our examples the sequence $x$ is a `transform'
${\T}(x_{\rm \scriptstyle F})$ of the Fibonacci word
$x_{\rm \scriptstyle F} = 01001010010010100\ldots$ These transforms have been introduced
by Kimberling in the OEIS \cite{oeis}. Our main example is:  $\T: [001\mapsto 2]$.
Replacing, in $x_{\rm \scriptstyle F} = 01001010010010100\ldots$, each $001$ by $2$ gives
$x_{\rm \scriptstyle K}=01201220120\ldots$.

\medskip

For the transform method $\T$  we can derive a `general' result similar to
Theorem~\ref{th:GBS-returns}. However, since Kimberling applies the StringReplace
procedure from Mathematica, which replaces occurrences of $w$ consecutively from left to
right, we do not obtain a sequence of return words in the case that $w$ has overlaps in
$x_{\rm \scriptstyle F}$. This restricts the number of words $w$ to which the following
theorem (Theorem~\ref{th:Kim}) applies. We first recall a definition from
\cite[p.~593--594]{WenWen}.

\begin{definition}[Wen \& Wen]
Let $w$ be a factor of the Fibonacci word $x_{\rm \scriptstyle F}$. We say that $w$ has an
\emph{overlap} in $x_{\rm \scriptstyle F}$ if there exist non-empty words $x, y$ and $z$ such
that $w = xy = yz$, and the word $xyz$ is a factor of the Fibonacci word.
\end{definition}

\begin{theorem}\label{th:Kim}
Let $x_{\rm \scriptstyle F}$ be the Fibonacci word, and let $w$ be a factor of
$x_{\rm \scriptstyle F}$ that has no overlap in $x_{\rm \scriptstyle F}$. Consider the
transform ${\T}(x_{\rm \scriptstyle F})$, which replaces every occurrence of the word
$w$ in $x_{\rm \scriptstyle F}$ by the letter $2$.
Let $Y$ be the sequence $({\T}(x_{\rm \scriptstyle F}))^{-1}(2)$, i.e., the positions of $2$'s
in ${\T}(x_{\rm \scriptstyle F})$. Then $Y$ is a generalized Beatty sequence, i.e., for all
$n\ge 1$, $Y(n) = p\lfloor n \varphi \rfloor + qn +r$, with parameters $p,q,r$, which can be
explicitly computed.
\end{theorem}

\proof As in the proof of Theorem~\ref{th:GBS-returns}, let
$x_{\rm \scriptstyle F}=r_0(w)r_1(w)r_2(w)\dots $, written as a concatenation of return
words of the word $w$. Now the distances between $2$'s in ${\T}(x_{\rm \scriptstyle F})$
are equal to $l_1:=|r_1(w)|-|w|+1$ and $l_2:=|r_2(w)|-|w|+1$. We can apply
Lemma~\ref{lem:diff}, which gives  $p=l_1-l_2, \; q=2l_2-l_1$. Inserting $n=1$, we find
that $r=|r_0(w)|-l_2+1$. \endpf

\begin{example}\label{mapsto2}
We take  $\T:[001\mapsto 2]$, with image ${\T}(x_{\rm \scriptstyle F})=
01201220120\dots$, so $Y=(3,6,7,10,\dots)$. Here $r_0(w)=01, r_1(w)=00101, r_2(w)=001$.
This gives $l_1=3,\, l_2=1$, implying $p=2$, $q=-1$ and $r=2$. So $Y$ is
the generalized Beatty sequence $(Y_n)_{n\ge 1} = (2\lfloor{n\varphi}\rfloor-n+2)_{n\ge 1}$.
\end{example}

The question arises whether not only ${\T}(x_{\rm \scriptstyle F})^{-1}(2)$, but also
${\T}(x_{\rm \scriptstyle F})^{-1}(0)$ and ${\T}(x_{\rm \scriptstyle F})^{-1}(1)$  are
generalized Beatty sequences. In general this is not true. However, this holds for
$\T:[001\mapsto 2]$.

\begin{theorem}\label{th:z-prime}
Let $\T : [001 \mapsto 2]$, and let $x_{\rm \scriptstyle K} := \T (x_{\rm \scriptstyle F})$.
Then the three sequences $x_{\rm \scriptstyle K}^{-1}(0)$, $x_{\rm \scriptstyle K}^{-1}(1)$,
$x_{\rm \scriptstyle K}^{-1}(2)$ form a complementary triple of generalized Beatty sequences.
\end{theorem}

\proof
According to Example~\ref{mapsto2} we have that
$x_{\rm \scriptstyle K}^{-1}(2) = (2 \lfloor n \varphi \rfloor - n + 2)_{n \geq 1}$.
Since clearly $x_{\rm \scriptstyle K}^{-1}(0) = x_{\rm \scriptstyle K}^{-1}(1) - 1$, our remaining
task is to prove that $(Z(n))_{n \geq 0} := x_{\rm \scriptstyle K}^{-1}(1)$ is a generalized
Beatty sequence. The return word structure of the word $w = 001$ in $x_{\rm \scriptstyle F}$
is given by
$$
r_0(w) = 01, \ \ r_1(w) = 00101, \ \ r_2(w) = 001.
$$
Note that $Z(0) = 2$, the $1$ coming from $r_0(w)$. This is exactly the reason why it is
convenient to start $Z$ from index $0$: the other $1$'s are coming from the $r_1(w)$'s---note
that $r_2(w)$ is mapped to $2$.

The differences between the indices of occurrences of $2$ in $x_{\rm \scriptstyle K}$ are given
by the Fibonacci word $3133131331\ldots$, which codes the appearance of the words $r_1(w)$
and $r_2(w)$. Therefore, to obtain the differences between the indices of occurrences of $1$ in
$x_{\rm \scriptstyle K}$, we have to map the word $w' = 13$ to $4$ , obtaining the word
$u = 343443\ldots$. To obtain a description of $u$, we apply Theorem~\ref{th:GBS-returns} a
second time with $w' = 13$. We have $r_0(w') = 3$, $r_1(w') = 133$, $r_2(w') = 13$.
So $ l_1 = |r_1(w')| - |w'| + 1 = 2$, and $l_2 = |r_2(w')| - |w'| + 1 = 1$, which
give  $p =l_1-l_2= 1$, $q = 2l_2-l_1=0$. The conclusion is that positions of $4$ in $u$ are given by the
generalized Beatty sequence $(\lfloor n \varphi \rfloor + 1)_{n \geq 1}$. This forces that $u$ is
nothing else than the Fibonacci word on $\{4, 3\}$, preceded by $3$. But then $Z$ is a
generalized Beatty sequence with parameters $p = 1$, $q = 2$. Since $Z(1) = 5$, we must have
$r = 2$, which happens to fit perfectly with the value $Z(0) = 2$. \endpf

\medskip

\noindent Here is an example where ${\T}(x_{\rm \scriptstyle F})^{-1}(0)$ and ${\T}(x_{\rm \scriptstyle F})^{-1}(1)$  are \emph{not} generalized Beatty sequences.

\begin{example}\label{ex:00100}
We take  $\T:[00100\mapsto 2]$, with image ${\T}(x_{\rm \scriptstyle F})=
010010121010010121012\dots$, so $Y=(8,17,21\dots)$.
Here $r_0(w)=0100101, r_1(w)=0010010100101, r_2(w)=00100101$.
This gives $l_1=9, l_2=4$, so $p=5$ and $q=-1$ and $r=4$. So $Y$ is the
generalized Beatty sequence $(Y_n)_{n\ge 1} = (5\lfloor{n\varphi}\rfloor-n+4)_{n\ge 1}$.
The positions of $0$ are given by
$({\T}(x_{\rm \scriptstyle F}))^{-1}(0)=1,3,4,6,10,12,13,\dots$, with difference sequence
$2,1,2,4,2,1,\dots$, so by Lemma~\ref{lem:diff} this sequence is not a generalized Beatty
sequence. However, it can be shown that $({\T}(x_{\rm \scriptstyle F}))^{-1}(0)$ is a union
of $4$ generalized Beatty sequences, and the same holds for
$({\T}(x_{\rm \scriptstyle F}))^{-1}(1)$.
\end{example}

%%OLD VERSION%%\noindent Example~\ref{ex:00100} raises the question whether the sequences
%%OLD VERSION%%${\T}(x_{\rm \scriptstyle F})^{-1}(0)$ and ${\T}(x_{\rm \scriptstyle F})^{-1}(1)$ are always
%%OLD VERSION%%finite unions of generalized Beatty sequences. This can be proved---generalizing the proof
%%OLD VERSION%%of Theorem~\ref{th:z-prime}---under the condition that
%%OLD VERSION%%$$|r_0(w)|\le |r_1(w)|-|w|\qquad {\rm (SR0)}.$$

\noindent Here is the general result.

\begin{theorem}\label{th:unionGBS}
For a non-overlapping word $w$ from the Fibonacci language let $\T : [w \mapsto 2]$, and let $Y := \T (x_{\rm \scriptstyle F})$.
Suppose $w$ satisfies
$$|r_0(w)|\le |r_1(w)|-|w|\qquad {\rm (SR0)}.$$
Then the three sequences $Y^{-1}(0)$, $Y^{-1}(1)$,  $Y^{-1}(2)$
are finite unions of generalized Beatty sequences.
\end{theorem}

\noindent Note that we already know by Theorem \ref{th:Kim} that  $Y^{-1}(2)$ is a single generalized Beatty sequence. The condition (SR0) states that the length of $r_0(w)$ is small with respect to the length of $r_1(w)$.

For the proof of Theorem \ref{th:unionGBS} one needs the following proposition.

\begin{proposition}\label{prop:r0r1r2} Let $w$ be a word from the Fibonacci language, and let
$r_0(w)r_1(w)r_2(w)\dots$ be the return sequence of $w$ in the Fibonacci word
$x_{\rm \scriptstyle F}$. Then {\rm (i)} $r_0(w)$ is a suffix of $r_1(w)$, and {\rm (ii)}
if $r_2(w)=wt_2(w)$, then $t_2(w)$ is a suffix of $r_1(w)$.
\end{proposition}

\proof  Let $s_0=1, s_1=00, s_2=101, s_3=00100,\dots$ be the singular words introduced
in \cite{WenWen}. According to  \cite[Theorem~1.9.]{HuangWen} there is a unique largest
singular word $s_k$ occurring in $w$, so we can write $w=\mu_1 s_k \mu_2$,
for two words $\mu_1, \mu_2$ from the Fibonacci language. It is known---see
\cite{WenWen}  and the remarks after \cite[Proposition~1.6.]{HuangWen}--- that the
two return words of the singular word $s_k$ are
 $$r_1(s_k)=s_ks_{k+1}, \qquad r_2(s_k)=s_ks_{k-1}.$$
 According to \cite[Lemma~3.1]{HuangWen}, the two return words of $w$ are given by
 $$r_1(w)=\mu_1 r_1(s_k) \mu_1^{-1},\quad  r_2(w)=\mu_1 r_2(s_k) \mu_1^{-1}.$$
 Substituting the first  equation in the second, we obtain the key equation
 \begin{equation}\label{eq:rw}
 r_1(w)=\mu_1 s_ks_{k+1} \mu_1^{-1},\quad  r_2(w)=\mu_1 s_ks_{k-1} \mu_1^{-1}.
 \end{equation}

\noindent Proof of (i): \quad
We compare the return word decompositions of $x_{\rm \scriptstyle F}$
by $s_k$ and by $w$:
$$r_0(s_k)r_1(s_k)r_2(s_k)r_1(s_k)\dots = r_0(w)r_1(w)r_2(w)r_1(w)\dots=
r_0(w)\mu_1 r_1(s_k) \mu_1^{-1}\mu_1 r_2(s_k)\mu_1^{-1}\mu_1 r_1(s_k)\mu_1^{-1}
\dots$$
 It follows that we must have $r_0(s_k)=r_0(w)\mu_1$, and so
 $r_0(w)=r_0(s_k)\mu_1^{-1}$.
By \cite[Lemma~2.3]{HuangWen}, $r_0(s_k)$ equals $s_{k+1}$, with the first letter deleted.
Thus we obtain from Equation~(\ref{eq:rw}) that $r_0(w)$ is a suffix of $r_1(w)$.

 \medskip

\noindent Proof of (ii): \quad
Since $s_{k+1}=s_{k-1}s_{k-3}s_{k-1}$, by  \cite[Property~2]{WenWen}, we can do the
following computation, starting from Equation~(\ref{eq:rw}):
 $$r_1(w)=\mu_1 s_ks_{k+1}\mu_1^{-1}=w\mu_2^{-1}s_{k+1}\mu_1^{-1}=
w\mu_2^{-1} s_{k-1}s_{k-3}s_{k-1}\mu_1^{-1}=
w\mu_2^{-1} s_{k-1}s_{k-3}\mu_2\mu_2^{-1}s_{k-1}\mu_1^{-1}. $$
For $r_2(w)$ we have
$$r_2(w)=\mu_1 s_ks_{k-1}\mu_1^{-1}=w\mu_2^{-1}s_{k-1}\mu_1^{-1}. $$
Now note that in this concatenation $\mu_2^{-1}$ cancels against a suffix of $w$.
We claim that it also cancels against a prefix of $s_{k-1}$. This follows, since by
\cite[Proposition 2.5]{HuangWen} \emph{any} occurrence of $s_k$ in
$x_{\rm \scriptstyle F}$ is directly followed by a $s_{k+1}=s_{k-1}s_{k-3}s_{k-1}$ with the
last letter deleted. It now follows that $t_2(w)=\mu_2^{-1}s_{k-1}\mu_1^{-1}$, and we see
that this word is a suffix of $r_1(w)$. \endpf

\medskip

\noindent {\emph{Proof of Theorem~\ref{th:unionGBS}:}}

\noindent From Property (ii) in Proposition \ref{prop:r0r1r2} we obtain that the return words of $w$ can be written as
$$ r_1(w) = w\, m_1(w)\, t_2(w), \quad  r_2(w) = w\, t_2(w),$$
for some words $ m_1(w)$ and $ t_2(w)$. Let $Z:=Y^{-1}(2)$ be the the positions of the letter  2.
If $t_2(w)$ is non-empty, then any letter in $t_2(w)$ occurs in $Y := \T (x_{\rm \scriptstyle F})$ in positions which are just a shift $ -|t_2(w)|,\dots,-1$ of $Z$, so each letter occurs according to a generalized Beatty sequence. The word $ m_1(w)$ is never empty, and any letter in $m_1(w)$ occurs in $Y = \T (x_{\rm \scriptstyle F})$ in positions which are a shift of a subsequence of $Z$ (except, possibly, for the first occurrence, which then is in $r_0(w))$. This subsequence is obtained by replacing the distances $\ell_1$ and $\ell_2$ of $\Delta Z$ by $\ell_1+\ell_2$ and $\ell_1$. Moreover, these distances occur as the Fibonacci word $x_{\rm \scriptstyle F}$ on the alphabet $\{\ell_1+\ell_2,\ell_1\}$, because $x_{\rm \scriptstyle F}$ is invariant under $0\mapsto 01, 1\mapsto 0$. Thus each letter in $m_1(w)$ occurs according to a generalized Beatty sequence. All these $|t_2(w)|+|m_1(w)|$ sequences start at index 1. If we let the last $|r_0(w)|$ of these sequences start at index 0, then we have taken into account \emph{all} elements of $Y$. This works, because of Property (i) in Proposition \ref{prop:r0r1r2}.   \endpf

\medskip

\noindent Here is an example where the (SR0) condition is not satisfied.

\smallskip

\begin{example}\label{ex:10100}
We take  $\T:[10100\mapsto 2]$, with image ${\T}(x_{\rm \scriptstyle F})=
01002100221002\dots$, so $Y=(5,9,10\dots)$. Here $r_0(w)=0100, \,r_1(w)=10100100,\,
r_2(w)=10100$.
%This gives $l_1=9, l_2=4$, so $p=5$ and $q=-1$ and $r=4$. So \phantom{}{$Y$} is
%the generalized Beatty sequence $(Y_n)_{n\ge 1} =
%(5\lfloor{n\varphi}\rfloor-n+4)_{n\ge 1}$.
The positions of $0$ are given by the sequence $({\T}(x_{\rm \scriptstyle F}))^{-1}(0)=1,3,4,7,8\dots$,
which can be written as a union of two  generalized Beatty sequences, \emph{except} that
the position $1$ from the first $0$ in ${\T}(x_{\rm \scriptstyle F})$ will not be in this union.
\end{example}

With Equation~(\ref{eq:rw}) we can deduce an equivalent simple formulation of condition
(SR0). If $w=\mu_1s_k\mu_2$, then $r_0(w)$ equals $s_{k+1}\mu_1^{-1}$ with the first
letter removed, and $r_1(w)=\mu_1 s_ks_{k+1} \mu_1^{-1}$, so
$$|w|=|\mu_1|+F_k+|\mu_2|,\quad |r_0(w)|=F_{k+1}-|\mu_1|-1,\quad |r_1(w)|=
F_{k+1}+F_k.$$
Filling this into condition (SR0) we obtain
$$|\mu_2|\le 1\qquad {\rm (SR0')}.$$
Using ${\rm (SR0')}$, together with  Theorem~6 in \cite{WenWen}, one can show that  Theorem~\ref{th:unionGBS} does apply to at most $3$ words $w$ of length
$m$, for all $m \geq 2$ (in fact, only $2$, if $m$ is not a Fibonacci number).

\bigskip

\noindent
{\bf Acknowledgments} We thank the referee for useful comments.

\end{document}